\documentclass[11pt,reqno]{amsart}
\usepackage{graphicx}
\usepackage{amsfonts}
\usepackage{mathrsfs}
 \usepackage{times}
 \usepackage{mathptmx}
\usepackage{amsmath}
\usepackage{txfonts}
\usepackage{amssymb}
\usepackage{amsmath,amsthm}
\usepackage{CJK,amsfonts}
\usepackage[body={7.5in,8.75in}, top=1.2in, right=0.9in,left=0.9in]{geometry}
\usepackage[colorlinks=true, linkcolor=red, citecolor=blue]{hyperref}
\allowdisplaybreaks[4] \numberwithin{equation}{section}

\newtheorem{thm}{Theorem}

\newtheorem{corollary}[thm]{Corollary}
\newtheorem{lemma}[thm]{Lemma}
\newtheorem{remark}[thm]{Remark}
\newtheorem{proposition}[thm]{Proposition}

\def\squarebox#1{\hbox to #1{\hfill\vbox to #1{\vfill}}}

\newcommand{\qbinomial}[3]{\mbox{$
\biggl[ 
\begin{array}{c}
#1\\
 #2
\end{array}\biggr]_{
\!{#3}}$}}

\newcommand{\be}{\begin{equation}}
\newcommand{\ee}{\end{equation}}
\newcommand{\bea}{\begin{eqnarray}}
\newcommand{\eea}{\end{eqnarray}}
\newcommand{\bd}{\begin{displaymath}}
\newcommand{\ed}{\end{displaymath}}
\begin{document}
\begin{CJK*}{GBK}{song}
\title[$q$-difference equation for    generalized trivariate $q$-Hahn   polynomials]{$q$-difference equation for    generalized trivariate $q$-Hahn   polynomials}
\author{  Sama Arjika${}^{1}$ and Mahaman Kabir Mahaman${}^2$}
\dedicatory{\textsc}
\thanks{${}^1$Department of Mathematics and Informatics, University of Agadez, Niger.  ${}^2$Department of Mathematics,
ENS,  University Abdou Moumouni of  Niamey, B.P.: 10 896, Niamey,  Niger}
\thanks{Emails:  rjksama2008@gmail.com (S.A.), mahamankabir@yahoo.fr (M.K.M)}



\begin{abstract}In this   paper, we introduce   a family of   trivariate $q$-Hahn   polynomials $\Psi_n^{(a)}(x,y,z|q)$ as a general form of Hahn polynomials $\psi_n^{(a)}(x|q),$ $\psi_n^{(a)}(x,y|q)$ and $F_n(x,y,z;q)$.   We represent $\Psi_n^{(a)}(x,y,z|q)$ by the homogeneous $q$-difference operator $\widetilde{L}(a,b; \theta_{xy})$ introduced  by   Srivastava {\it et al} [H.  M. Srivastava,  S. Arjika  and A.  Sherif Kelil,   {\it  Some homogeneous $q$-difference operators and the associated  generalized Hahn polynomials}, Appl. Set-Valued Anal. Optim. {\bf 1} (2019), pp. 187--201.] to derive: extended generating, Rogers formula, extended Rogers formula and Srivastava-Agarwal type generating functions involving $\Psi_n^{(a)}(x,y,z|q)$  by the $q$-difference equation.  
\end{abstract}

\maketitle

\section{\bf Introduction} 

In this paper, we adopt the common conventions and notations on $q$-series.  For the convenience of the reader, we provide  a  summary of the mathematical notations, basics properties  and definitions to be used  in the sequel.   We refer  to the general references (see \cite{Koekock})  for the definitions and notations. Throughout this paper, we assume that  $|q|< 1$. 

For complex  numbers $a$, the $q$-shifted factorials are defined by:
\be 
(a;q)_0:=1,\quad  (a;q)_{n}: =\prod_{k=0}^{n-1} (1-aq^k),   \quad (a;q)_{\infty}:=\prod_{k=0}^{\infty}(1-aq^{k})
\ee 
and 
 $ (a_1,a_2, \ldots, a_r;q)_m=(a_1;q)_m (a_2;q)_m\cdots(a_r;q)_m,\; m\in\{0, 1, 2\cdots\}$.  
 \\
 The $q$-binomial coefficient is defined as \cite{GasparRahman}
$$
 {\,n\,\atopwithdelims []\,k\,}_{q}=\frac{(q;q)_n}{(q;q)_k(q;q)_{n-k}}=\frac{(q^{-n};q)_k}{(q;q)_k}(-1)^kq^{nk-({}^k_2)},\,\, \mbox{ for } 0\leq k\leq n.
$$
The    basic or $q$-hypergeometric function 
  in the variable $z$ (see  Slater \cite[Chap. 3]{SLATER},  Srivastava and Karlsson   \cite[p. 347, Eq. (272)]{SrivastaKarlsson}   for details) is defined as:
 $$
{}_{r}\Phi_s\left[\begin{array}{c}a_1, a_2,\ldots, a_r;
 \\\\
b_1,b_2,\ldots,b_s;
 \end{array}
q;z\right]
 =\sum_{n=0}^\infty\Big[(-1)^n q^{({}^n_2)}\Big]^{1+s-r}\,\frac{(a_1, a_2,\ldots, a_r;q)_n}{(b_1,b_2,\ldots,b_s;q)_n}\frac{ z^n}{(q;q)_n},
$$
 when $r>s+1$. Note that, for $r=s+1$, we have:
$$
{}_{r+1}\Phi_r\left[\begin{array}{r}a_1, a_2,\ldots, a_{r+1};
 \\\\
b_1,b_2,\ldots,b_r;
 \end{array}q;z\right]
 =\sum_{n=0}^\infty \frac{(a_1, a_2,\ldots, a_{r+1};q)_n}{(b_1,b_2,\ldots,b_r;q)_n}\frac{ z^n}{(q;q)_n}.
$$

We will be mainly concerned with the Cauchy polynomials as given below  \cite{Chen2003}
\bea
\label{def}
 p_n(x,y):=(x-y)( x- qy)\cdots ( x-q^{n-1}y) =(y/x;q)_n\,x^n
\eea
with the  Srivastava-Agarwal type generating function  
\be
\label{Srivas}
\sum_{n=0}^\infty  
p_n (x,y)\frac{(\lambda;q)_n\,t^n}{(q;q)_n}= {}_{2}\Phi_1\left[\begin{array}{r}\lambda,y/x;
 \\\\
0;
 \end{array}
q; xt\right].
 \ee
 For $\lambda=0,$ we get the generating function \cite{Chen2003}
\be
\label{gener}
\sum_{n=0}^{\infty} p_n(x,y)
\frac{t^n }{(q;q)_n} = 
\frac{(yt;q)_\infty}{(xt;q)_\infty}.
\ee
 The generating   function (\ref{gener}) is also the   homogeneous version  of the Cauchy identity or the $q$-binomial theorem   given by \cite{GasparRahman}
\be
\label{putt}
\sum_{k=0}^{\infty} 
\frac{(a;q)_k }{(q;q)_k}z^{k}={}_{1}\Phi_0\left[\begin{array}{r}a;\\
 \\
-;
 \end{array}
q;z\right]= 
\frac{(az;q)_\infty}{(z;q)_\infty},\quad |z|<1. 
\ee
Putting    $a=0$, the relation (\ref{putt}) becomes  Euler's identity  \cite{GasparRahman}
\be
\label{q-expo-alpha}
  \sum_{k=0}^{\infty} \frac{ z^{k}}{(q;q)_k}=\frac{1}{(z;q)_\infty}\quad |z|<1
\ee
and its inverse relation  \cite{GasparRahman}
\be
\label{q-Expo-alpha}
 \sum_{k=0}^{\infty}  \frac{(-1)^kq^{ ({}^k_2)
}\,z^{k}}{(q;q)_k}=(z;q)_\infty.
\ee
 
    Saad and  Sukhi \cite{Saadsukhi} defined the   $q$-difference operator ${\theta}_{xy}$ 
\be 
\label{deffd}
  {\theta}_{xy}\big\{f(x,y)\}:=\frac{f(q^{-1}x,y)-f( x,qy)}{q^{-1}x-y},
\ee
which  turns out to be suitable for dealing with the Cauchy polynomials.  Their corresponding $q$-exponential operator is
\be 
\label{qoperator}
\mathbb{E}(z\theta_{xy})=\sum_{k=0}^\infty \frac{q^{({}^k_2)} }{(q;q)_k}\, \left(z\,\theta_{xy}\right)^k.
\ee

Recently, Srivastava, Arjika and Kelil \cite{HariSama} have introduced  the $q$-difference operator  $\widetilde{L}(a,b; \theta_{xy})$   
\be 
\label{operator}
\widetilde{L}(a,b; \theta_{xy})=\sum_{k=0}^\infty \frac{q^{({}^k_2)}\,(a;q)_k }{(q;q)_k}\, \left(b\,\theta_{xy}\right)^k,
\ee
to study $q$-polynomials and related generating functions. 

In this paper, our goal is to generalize the results of Srivastava, Arjika and Kelil \cite{HariSama}, and  Mohameed  \cite{Mohammed}. We first construct the following  generalized trivariate $q$-Hahn   polynomials  as
\be 
\label{ndef}
\Psi_n^{(a)}(x,y,z|q)=(-1)^nq^{-({}^n_2)} \sum_{k=0}^n {\,n\,\atopwithdelims []\,k\,}_{q}(-1)^kq^{({}^k_2)} (a ;q)_k  p_{n-k}(y,x)z^k. 
\ee
\begin{remark}
 For $a=0 $,  the  generalized trivariate $q$-Hahn polynomials $\Psi_n^{(a)}(x,y,z|q)$ are the well known trivariate $q$-polynomials $F_n(x,y,z;q)$ investigated by Mohameed (see \cite{Mohammed} for more details), i.e.,
\be
\Psi_n^{(0)}(x,y,z|q) = F_n(x,y,z;q).
\ee
If we let $a=0,\, y = ax$ and $z=y$, the  generalized trivariate $q$-Hahn   polynomials $\Psi_n^{(a)}(x,y,z|q)$ reduce to  the second Hahn polynomials $\psi_n^{(a)}(x,y|q)$ \cite{CaoJ12}, i.e.,
\be
\Psi_n^{(0)}(x,ax,y|q) =   \psi_n^{(a)}(x,y|q).
\ee
Also, $a=0,\, y = ax$ and $z=1$,  the  generalized trivariate $q$-Hahn $\Psi_n^{(a)}(x,y,z|q)$   reduce to   Hahn polynomials $\psi_n^{(a)}(x|q)$ \cite{AlSalam}, i.e., 
\be
\Psi_n^{(0)}(x,ax,1|q) = \psi_n^{(a)}(x|q).
\ee
\end{remark}
The polynomials (\ref{ndef}) can be represented by the homogeneous $q$-difference 
operator (\ref{operator}) as follows.
\begin{proposition}
\be
\label{conds}
\Psi_n^{(a)}(x,y,z|q)=\widetilde{L}(a,z; \theta_{xy})\left\{(-1)^nq^{-({}^n_2)} p_n(y,x)\right\}.
\ee
\end{proposition}
 \begin{proof}
 By identity (\ref{operator}) and  taking into account $ \displaystyle  \theta_{xy}p_n(y,x)=-(1-q^n)\,p_{n-1}(y,x)$,  we get the result.  
 \end{proof}
 In light of 
 $ \displaystyle 
    \theta_{xy}^k [(xt;q)_\infty/(yt;q)_\infty]=(-t)^k[ (xt;q)_\infty/(yt;q)_\infty],$
we have the following identity
\be
\label{identity}
\widetilde{L}(a,z; \theta_{xy})\left\{\frac{(xt;q)_\infty}{(yt;q)_\infty}\right\}=\frac{(xt;q)_\infty}{(yt;q)_\infty}\;{}_{1}\Phi_1\left[\begin{array}{r}a; 
 \\\\
0 ;
 \end{array} 
q;zt\right].
\ee

The main object of this paper is  to  use the $q$-difference equation  to  derive some identities such as: extended generating function, Rogers  formula, extended Rogers   formula and   Srivastava-Agarwal type generating functions. 
 
The  paper is organized as follows. In Section \ref{qdifff}, we   state two theorems and give the proofs.  We derive  an  extended generating function  for  these $q$-polynomials.    In Section \ref{section3}, we state the  Rogers  formula and extended Rogers  formula and give the proofs by the  $q$-difference equation.   In  Section \ref{qdiffff}, we obtain  Srivastava-Agarwal type generating functions involving the  generalized trivariate $q$-Hahn  polynomials by the method of $q$-difference equation.

\section{Main results and proofs}
\label{qdifff}
In this section, we introduce another extension of  $q$-Hahn   polynomials. Then, we represent it by the homogeneous $q$-difference  operator and derive  an extended generating function. 
\begin{thm}
\label{aprodpos}
Let $f (a, b,x,y, z)$ be an $5$-variable analytic function at 
  $(a, b,x,y, z)=(0,0,0,0, 0)\in\mathbb{C}^{5}$. If $f (a, b,x,y, z)$ satisfies the $q$-difference equation
\begin{multline}
\label{asaaaa}
 (q^{-1}x- y)\Big[f(a, b,x,y, z)-f(a,  b,x,y,qz)\Big]=  z \Big[f(a,  b,q^{-1}x,y, qz)-\\
 f(a,  b,x, qy, qz)\Big]+ az \Big[f(a,  b,x, qy,q^2z)- f(a,   b,q^{-1}x,y,q^2z)\Big],
\end{multline}
then we have:
\bea
\label{shha}
f(a,  b,x,y, z)=\widetilde{L}(a,z; \theta_{xy})\Big\{f(a , b,x,y, 0) \Big\}.
\eea
\end{thm}
\begin{corollary}
\label{aprodpos}
Let $f ( b,x,y, z)$ be an $4$-variable analytic function at 
  $(b,x,y, z)=(0,0,0, 0)\in\mathbb{C}^{4}$. If $f (b,x,y, z)$ satisfies the $q$-difference equation
\begin{align}
\label{asaaaa}
 (q^{-1}x- y)\Big[f(b,x,y, z)-f(b,x,y,qz)\Big]
 =  z \Big[f(b,q^{-1}x,y, qz)-f(b,x, qy, qz)\Big],
\end{align}
then we have:
\bea
\label{shha}
f(b,x,y, z)=\mathbb{E}(z \theta_{xy})\Big\{f(b,x,y, 0) \Big\}.
\eea
\end{corollary}
\begin{proof}
 From the theory of several complex variables \cite{Range}, we begin to solve the $q$-difference equation  (\ref{asaaaa}). First we may assume that
\be
\label{aa120}
f(a,  b,x,y, z)= \sum_{ n=0}^\infty A_n(a, b,x,y)z^n. 
\ee
Substituting (\ref{aa120})  into (\ref{asaaaa}), we get:
\begin{align*}
(q^{-1}x- y)\sum_{n=0}^\infty (1-q^n)A_n(a, b, x,y)z^n 
=  \sum_{n=0}^\infty q^n(1-aq^n)\Big[A_n(a , b,q^{-1}x,y)
- A_n(a, b, x,qy)\Big]z^{n+1}. 
\end{align*}
 Comparing coefficients of $z^n,\,n\geq 1$, we   find that
\begin{align*}
(q^{-1}x- y)  (1-q^n)A_n(a, b,x,y) 
=  q^{n-1}(1-aq^{n-1})\Big[A_{n-1}(a,  b,q^{-1}x,y)-A_{n-1}(a, b, x,qy) 
\Big]. 
\end{align*}
After simplification, we get:
\be 
 A_n(a, b,x,y)=    q^{n-1}\frac{1 -aq^{n-1} }{1-q^n} \theta_{xy}\Big\{ A_{n-1}(a, b,x,y)\Big\}.\nonumber
\ee
By iteration, we gain
\bea
\label{aaZddA}
   A_n(a, b,x,y) 
 =   q^{({}^n_2)}  \frac{(a;q)_n }{(q;q)_n} \theta_{xy}^n \Big[A_0(a, b,x,y)\Big].
\eea
  Just taking $z=0$ in (\ref{aa120}), we immediately obtain
$A_0(a, b,x,y) =f(a, b,x,y,0)$. Substituting  (\ref{aaZddA}) back into (\ref{aa120}), we achieve  (\ref{asaaaa}). 
\end{proof}

\begin{thm}[Extended generating function for $\Psi_n^{(a)}(x,y,z|q)$] 
\label{dgend} For $|yt|<1,$ we have:
\begin{align}
\label{exdddgen}
\sum_{n=0}^\infty  \Psi_{n+k}^{(a)}(x,y,z|q)\frac{ (-1)^{n+k}q^{ ({}^{n+k}_{\,\,\,2})}\,t^n}{(q;q)_n}
= t^{-k}\frac{ (xt;q)_\infty}{ (yt;q)_\infty}\sum_{n=0}^k\frac{ (q^{-k},yt;q)_n\,q^n}{(xt,q;q)_n} \;{}_{1}\Phi_1\left[\begin{array}{r}a ;
 \\\\
0;
 \end{array}
q;  ztq^n\right].
\end{align}
\end{thm}
\begin{corollary} For $|yt|<1,$ we have:
\begin{align}
\label{dgen}
\sum_{n=0}^\infty  F_{n+k}(x,y,z;q)\frac{ (-1)^{n+k}q^{ ({}^{n+k}_{\,\,\,2})}\,t^n}{(q;q)_n}
= t^{-k}\frac{ (xt,zt;q)_\infty}{ (yt;q)_\infty}\;{}_{3}\Phi_2\left[\begin{array}{r}q^{-k},yt,0 ;
 \\\\
xt,zt;
 \end{array}
q;  q\right].
\end{align}
\end{corollary}
\begin{remark}
For $a=0$, (\ref{exdddgen}) reduces (\ref{dgen}). 
For $a=0$ and $k=0$, (\ref{exdddgen}) and (\ref{dgen}) reduce  to  the generating function for $\Psi_n^{(a)}(x,y,z|q)$ 
\be 
\label{gen}
\sum_{n=0}^\infty  \Psi_n^{(a)}(x,y,z|q)\frac{ (-1)^nq^{ ({}^n_2)}\,t^n}{(q;q)_n}= \frac{(xt;q)_\infty}{(yt;q)_\infty}\;{}_{1}\Phi_1\left[\begin{array}{r}a; 
 \\\\
0;
 \end{array} 
q;  zt\right],\quad |yt|<1 
\ee 
and \cite[Theorem 2.6]{Mohammed}.
\end{remark}
To prove   the  Theorem \ref{dgend}, the  following Lemma is necessary. 
\begin{lemma}
$q$-Chu-Vandermonde  formula \cite[Eq. (II.6)]{GasparRahman}
\be
\label{qchuv}
{}_{2}\Phi_1\left[\begin{array}{r} q^{-n}, a;
 \\\\
c;
 \end{array} 
q; q \right]=\frac{(c/a;q)_n}{(c;q)_n}a^n.
\ee
\end{lemma}
\begin{proof}[Proof of Theorem \ref{dgend}]  
Denoting the right-hand side of equation (\ref{exdddgen}) by $F(a,t,x,y,z)$, we have:
\be 
\label{tass}
F(a,t,x,y,z)=t^{-k}\sum_{n=0}^k\frac{ (q^{-k};q)_n\,q^n}{(q;q)_n}\frac{ (xtq^n;q)_\infty}{ (ytq^n;q)_\infty} \;{}_{1}\Phi_1\left[\begin{array}{r}a ;
 \\\\
0;
 \end{array} 
q;  ztq^n\right].
\ee 
Because equation (\ref{tass}) satisfies (\ref{asaaaa}),   we have:
\bea
\label{ss}
F(a,t,x,y,z)&=&\widetilde{L}(a,z; \theta_{xy})   \left\{
F(a,t,x,y,0)\right\} \cr
&=&\widetilde{L}(a,z; \theta_{xy})\Bigg\{ t^{-k}\sum_{n=0}^k\frac{ (q^{-k};q)_n\,q^n}{(q;q)_n}\frac{ (xtq^n;q)_\infty}{ (ytq^n;q)_\infty}\Bigg\}
\cr
&=& \widetilde{L}(a,z; \theta_{xy})\Bigg\{\frac{ (xt;q)_\infty}{ (yt;q)_\infty}  t^{-k}\sum_{n=0}^k\frac{ (q^{-k},yt;q)_n\,q^n}{(xt,q;q)_n} \Bigg\}  
\cr
&=& \widetilde{L}(a,z; \theta_{xy})\left\{\frac{ (xt;q)_\infty}{ (yt;q)_\infty}  t^{-k}{}_{2}\Phi_1\left[\begin{array}{r}q^{-k},yt ;
 \\\\
xt;
 \end{array}
q; q\right] \right\} \, \mbox{by (\ref{qchuv})}\cr
&=& \widetilde{L}(a,z; \theta_{xy})\left\{\frac{p_{ k}(y,x)}{(xt;q)_k}\frac{(xt;q)_\infty}{(yt;q)_\infty}\right\}\label{AAA} \,\mbox{by (\ref{gener})}\cr
&=&\widetilde{L}(a,z; \theta_{xy})\Bigg\{p_k (y,x)\sum_{n=0}^\infty  
\frac{p_n(y,xq^k) \,t^n}{(q;q)_n}\Bigg\} 
 \cr
 &=&\widetilde{L}(a,z; \theta_{xy})\left\{\sum_{n=0}^\infty  p_{n+k}(y,x)\frac{ t^n}{(q;q)_n}\right\}\cr
 &=&\sum_{n=0}^\infty \widetilde{L}(a,z; \theta_{xy})\left\{ (-1)^{n+k}q^{ -({}^{n+k}_{\,\,\,2})}p_{n+k}(y,x)\right\}  (-1)^{n+k}q^{ ({}^{n+k}_{\,\,\,2})}\,\frac{t^n}{(q;q)_n} \,\,\,\mbox{by (\ref{conds})}\cr
 &=&\sum_{n=0}^\infty   \Psi_{n+k}^{(a)}(x,y,z|q)\frac{ (-1)^{n+k}q^{ ({}^{n+k}_{\,\,\,2})}\,t^n}{(q;q)_n},\nonumber
\eea 
which is the left-hand side of  (\ref{exdddgen}). 
\end{proof}
\section{The Rogers formula for $\Psi_{n}^{(a)}(x,y,z|q)$}
\label{section3}
In this section, we give  Rogers  formula and extended Rogers formula for the generalized trivariate $q$-Hahn polynomials $\Psi_{n}^{(a)}(x,y,z|q)$ by using the homogeneous $q$-difference equations.  
\begin{thm}[Roger's-type formula for $\Psi_{n}^{(a)}(x,y,z|q)$]  We have:
\label{exgend}
\begin{align}
\label{exgen}
\sum_{n=0}^\infty\sum_{m=0}^\infty  \Psi_{n+m}^{(a)}(x,y,z|q)(-1)^{n+m}q^{ ({}^{n+m}_{\,\,\,2})}\frac{ t^n }{(q;q)_n }\frac{  s^m}{ (q;q)_m}
=\frac{(xs;q)_\infty}{(t/s,ys;q)_\infty} \sum_{k=0}^\infty\frac{ (ys;q)_k\,q^k}{(sq/t,xs,q;q)_k} \;{}_{1}\Phi_1\left[\begin{array}{r}a;
 \\\\
0;
 \end{array} 
q;  zsq^k\right],
\end{align}
where $max\{|t/s|,|ys|\}<1$.
\end{thm}
\begin{corollary}  We have:
\label{dxgend}
\begin{align}
\label{fxgen}
\sum_{n=0}^\infty\sum_{m=0}^\infty  F_{n+m}(x,y,z;q)(-1)^{n+m}q^{ ({}^{n+m}_{\,\,\,2})}\frac{ t^n }{(q;q)_n }\frac{  s^m}{ (q;q)_m}
=\frac{(xs,zs;q)_\infty}{(t/s,ys;q)_\infty}  \;{}_{4}\Phi_3\left[\begin{array}{r}yt,0,0,0 ;
 \\\\
sq/t,xs,zs;
 \end{array}
q;  q\right],
\end{align}
where $max\{|t/s|,|ys|\}<1$.
\end{corollary}
\begin{proof}Denoting the right-hand side of equation (\ref{exgen}) by $G(a,s,x,y,z)$, we have:
\begin{align}
\label{ttass}
G(a,s,x,y,z)=\frac{1}{(t/s;q)_\infty} \sum_{k=0}^\infty\frac{ q^k}{(sq/t,q;q)_k}\frac{ (xsq^k;q)_\infty}{ (ysq^k;q)_\infty} \;{}_{1}\Phi_1\left[\begin{array}{r}a; 
 \\\\
0;
 \end{array} 
q; zsq^k\right].
\end{align}
Because equation (\ref{ttass}) satisfies (\ref{asaaaa}), by (\ref{shha}), we have:
\bea
\label{ss}
G(a,s,x,y,z)&=&\widetilde{L}(a,z; \theta_{xy})   \left\{
G(a,s,x,y,0)\right\} \cr
&=&\widetilde{L}(a,z; \theta_{xy})\Bigg\{\frac{1}{(t/s;q)_\infty} \sum_{k=0}^\infty\frac{ q^k}{(sq/t,q;q)_k}\frac{ (xsq^k;q)_\infty}{ (ysq^k;q)_\infty}\Bigg\}
\cr
&=& \widetilde{L}(a,z; \theta_{xy})\Bigg\{\frac{(xs;q)_\infty}{(ys;q)_\infty} \sum_{k=0}^\infty\frac{( ys;q)_kq^k}{(xs,q;q)_k} \frac{1}{(t/s;q)_\infty(sq/t;q)_k}\Bigg\}  
\cr
&=& \widetilde{L}(a,z; \theta_{xy})\Bigg\{\frac{(xs;q)_\infty}{(ys;q)_\infty} \sum_{k=0}^\infty\frac{( ys;q)_kq^k}{(xs,q;q)_k}\frac{(-t/s)^k q^{-({}^k_2)-k} }{(tq^{-k}/s;q)_\infty} \Bigg\} 
\cr
&=& \widetilde{L}(a,z; \theta_{xy})\Bigg\{\frac{(xs;q)_\infty}{(ys;q)_\infty} \sum_{k=0}^\infty\frac{( ys;q)_kq^k(-t/s)^k q^{-({}^k_2)-k}}{(xs,q;q)_k} \sum_{n=0}^\infty \frac{(t/s)^nq^{-nk}}{(q;q)_{n}} \Bigg\}\cr
&=& \widetilde{L}(a,z; \theta_{xy})\Bigg\{\frac{(xs;q)_\infty}{(ys;q)_\infty} \sum_{k=0}^\infty\frac{( ys;q)_kq^k}{(xs,q;q)_k}\sum_{n=k}^\infty \frac{(-t/s)^n q^{({}^k_2)-nk}}{(q;q)_{n-k}}
 \Bigg\}\cr
 &=& \widetilde{L}(a,z; \theta_{xy})\Bigg\{ \frac{(xs;q)_\infty}{(ys;q)_\infty}\sum_{n=0}^\infty \frac{(t/s)^n}{(q;q)_n}  \sum_{k=0}^n\frac{ (q^{-n},ys;q)_k\,q^k}{(xs,q;q)_k}
 \Bigg\}  
 \cr
 &=& \widetilde{L}(a,z; \theta_{xy})\left\{ \frac{(xs;q)_\infty}{(ys;q)_\infty}\sum_{n=0}^\infty \frac{(t/s)^n}{(q;q)_n}   
\;{}_{2}\Phi_1\left[\begin{array}{r}q^{-n},ys;
 \\\\
xs;
 \end{array} 
q;  q\right] \right\}
 \cr
 &=&\sum_{n=0}^\infty \frac{t^n}{(q;q)_n}    \widetilde{L}(a,z; \theta_{xy})\left\{ \frac{(xs;q)_\infty}{(ys;q)_\infty}
s^{-n}\;{}_{2}\Phi_1\left[\begin{array}{r}q^{-n},ys;
 \\\\
xs;
 \end{array} 
q; q\right] \right\}  \,\mbox{ by (\ref{qchuv})}
\cr
 &=&\sum_{n=0}^\infty \frac{t^n}{(q;q)_n}    \widetilde{L}(a,z; \theta_{xy}) \left\{\frac{p_n(y,x)\,(xs;q)_\infty}{(xs;q)_n(ys;q)_\infty}\right\} \,\mbox{by (\ref{AAA})}  
 \cr
&=&\sum_{n=0}^\infty\frac{t^n}{(q;q)_n}\sum_{m=0}^\infty  \Psi_{n+m}^{(a)}(x,y,z|q)(-1)^{n+m}q^{ ({}^{n+m}_{\,\,\,2})} \frac{ \,s^m}{ (q;q)_m},\nonumber
\eea
which is the left-hand side of  (\ref{exgen}). 
\end{proof}
\begin{thm}[Extended Roger's-type formula for $\Psi_{n}^{(a)}(x,y,z|q)$] 
\label{Tsxt} We have:
\begin{multline}
\label{gextend}
\sum_{n=0}^\infty \sum_{m=0}^\infty \sum_{k=0}^\infty\Psi_{n+m+k}^{(a)}(x,y,z|q) \frac{(-1)^{n+m+k}q^{ ({}^{n+m+k}_{\,\,\,\,\,\,2})}\,t^n\,s^m\omega^k}{(q;q)_{n+m}(q;q)_{m} (q;q)_k}\\
= \frac{ (x\omega;q)_\infty}{ (s/t,t/\omega,y\omega;q)_\infty}   \sum_{j=0}^\infty\frac{  (y\omega;q)_j\, q^j }{(x\omega,q\omega/t,q;q)_j}    \;{}_{1}\Phi_1\left[\begin{array}{r}a;
 \\\\
0;
 \end{array}
q;  z\omega q^j\right],
\end{multline}
where  $max\{|s/t|,|t/\omega|,|y\omega|\}<1$.
\end{thm}
\begin{remark}
For $s=0$, (\ref{gextend}) reduces to (\ref{exgen}).
\end{remark}
\begin{proof}[Proof of Theorem \ref{Tsxt}] 
Denoting the right-hand side of equation (\ref{gextend}) by $H(a,\omega,x,y,z)$, we have:
\begin{align}
\label{tttass}
H(a,\omega,x,y,z)= \frac{ 1 }{ (s/t,t/\omega;q)_\infty}   \sum_{k=0}^\infty\frac{  q^k }{(q\omega/t,q;q)_k} \cdot\frac{ (x\omega q^k;q)_\infty}{ ( y\omega q^k;q)_\infty}    \;{}_{1}\Phi_1\left[\begin{array}{r}a ;
 \\\\
0;
 \end{array} 
q; z\omega q^k\right].
\end{align}
Because equation (\ref{tttass}) satisfies (\ref{asaaaa}), by (\ref{shha}), we have:
\bea
\label{fass}
H(a,\omega,x,y,z)&=&\widetilde{L}(a,z; \theta_{xy})   \left\{
H(a,\omega,x,y,0)\right\}\cr
&=&\widetilde{L}(a,z; \theta_{xy})\Bigg\{\frac{1}{ (s/t,t/\omega;q)_\infty}   \sum_{k=0}^\infty\frac{  \, q^k }{( q\omega/t,q;q)_k}  \frac{ (x\omega q^k;q)_\infty}{ ( y\omega q^k;q)_\infty}  \Bigg\}\cr
&=&\frac{1}{ (s/t;q)_\infty}\widetilde{L}(a,z; \theta_{xy})\Bigg\{\frac{1}{ (t/\omega;q)_\infty}   \sum_{k=0}^\infty\frac{  \, q^k }{( q\omega/t,q;q)_k}  \frac{ (x\omega q^k;q)_\infty}{ ( y\omega q^k;q)_\infty}  \Bigg\}\cr
&=&\frac{1}{ (s/t;q)_\infty}\widetilde{L}(a,z; \theta_{xy})\left\{F(a,\omega,x,y,0) \right\} \,\mbox{ by (\ref{qchuv})}\cr
&=&\frac{1}{ (s/t;q)_\infty} \sum_{k=0}^\infty\sum_{m=0}^\infty  \Psi_{m+k}^{(a)}(x,y,z|q)    \frac{(-1)^{m+k}q^{ ({}^{m+k}_{\,\,\,2})}t^m}{(q;q)_m} \frac{\omega^k}{(q;q)_k}\cr
&=& \sum_{n=0}^\infty \frac{(s/t)^n}{(q;q)_n}\sum_{k=0}^\infty\sum_{m=0}^\infty  \Psi_{m+k}^{(a)}(x,y,z|q)   \frac{(-1)^{m+k}q^{ ({}^{m+k}_{\,\,\,2})} t^m}{(q;q)_m} \frac{\omega^k}{(q;q)_k}\cr
&=& \sum_{k=0}^\infty\sum_{m=0}^\infty   \sum_{m=n}^\infty\Psi_{n+m}^{(a)}(x,y,z|q)   \frac{(-1)^{m+k}q^{ ({}^{m+k}_{\,\,\,2})}t^{m-n}}{(q;q)_m} \frac{s^n}{(q;q)_n} \frac{\omega^k}{(q;q)_k}. 
\eea
Replacing $m$ by $m+n$ in (\ref{fass}), we obtain:
\begin{align}
H(a,\omega,x,y,z)=\sum_{n=0}^\infty \sum_{m=0}^\infty \sum_{k=0}^\infty\Psi_{n+m+k}^{(a)}(x,y,z|q)\,\frac{ (-1)^{n+m+k}q^{ ({}^{n+m+k}_{\,\,\,\,\,\,2})}t^n\,s^m\omega^k}{(q;q)_{m}(q;q)_{n+m}(q;q)_k},\nonumber 
\end{align}
which is the left-hand side of  (\ref{gextend}). 
\end{proof}

\section{Srivastava-Agarwal type generating functions for generalized trivariate $q$-Hahn  polynomials}
\label{qdiffff}

The Hahn polynomials \cite{Hahn049,Hahn49}
 (or Al-Salam and Carlitz polynomials \cite{AlSalam}) are defined as 
\be 
\Phi_n^{(a)}(x|q)=\sum_{k=0}^n \qbinomial{n}{k}{q} (a;q)_k x^k.
\ee 
Srivastava and Agarwal gave the  following generating function.
  \begin{lemma}\cite[Eq. (3.20)]{SrivastavaAgarwal}
Suppose   that  $max\{|t|, |xt|\}<1$, we have:  
\be
\label{c1sums}
\sum_{n=0}^\infty \Phi_n^{(\alpha)}(x|q) \frac{(\lambda;q)_n  t^n}{(q;q)_n}= \frac{(\lambda t; q)_\infty }{(t;q)_\infty}   {}_2\Phi_1\left[
\begin{array}{r} \lambda, \alpha;\\\\
 \lambda  t; \end{array}q;xt   
\right].
\ee
 \end{lemma}
 The generating function (\ref{c1sums}) is called Srivastava-Agarwal type generating functions for the Al-Salam-Carlitz polynomials \cite{SrivastavaAgarwal}.  

In this section, we give Srivastava-Agarwal type generating function  for the generalized  trivariate $q$-Hahn polynomials  $ \Psi_n^{(a)}(x,y,z|q)$ by the  homogeneous $q$-difference equation. 
\begin{thm}
\label{TA1}
For   $  |y\nu t|<1$, we have:  
\begin{align}  
\label{1sums}
\sum_{n=0}^\infty \Psi_n^{(a)}(x,y,z|q)p_n(\nu,\mu) \frac{(-1)^nq^{({}^n_2)}\,t^n}{(q;q)_n} 
= \frac{( \mu/\nu, x\nu t; q)_\infty }{(y\nu t;q)_\infty}\sum_{n=0}^\infty  \frac{(y\nu t;q)_n\,(\mu/\nu)^n   }{(x\nu t, q;q)_n}\, \;{}_{1}\Phi_1\left[\begin{array}{r}a ;
 \\\\
0;
 \end{array}
q; z\nu tq^n\right].
\end{align} 
 \end{thm}
\begin{corollary}
\label{TA1}
For   $  |y\nu t|<1$, we have:  
\begin{align}  
\label{1sums}
\sum_{n=0}^\infty F_n(x,y,z;q)p_n(\nu,\mu) \frac{(-1)^nq^{({}^n_2)}\,t^n}{(q;q)_n} 
= \frac{( \mu/\nu, x\nu t,z\nu t; q)_\infty }{(y\nu t;q)_\infty}\, \;{}_{3}\Phi_2\left[\begin{array}{r}y\nu t,0,0 ;
 \\\\
x\nu t,z\nu t;
 \end{array}
q; \frac{\mu}{\nu}\right].
\end{align} 
 \end{corollary}
\begin{corollary}
\label{TA1c}
For   $ |yt|<1$, we have:  
\begin{align} 
\label{1csums}
\sum_{n=0}^\infty \Psi_n^{(a)}(x,y,z|q) (\lambda;q)_n \frac{(-1)^nq^{({}^n_2)}\,t^n}{(q;q)_n} 
= \frac{(\lambda, x  t; q)_\infty }{(y  t;q)_\infty}\sum_{n=0}^\infty  \frac{(y  t;q)_n\,\lambda^n   }{(x  t, q;q)_n}\, \;{}_{1}\Phi_1\left[\begin{array}{r}a; 
 \\\\
0;
 \end{array} 
q; z  tq^n\right].
\end{align}
 \end{corollary}
 \begin{remark}
 For $\nu=1$, (\ref{1sums}) reduces to (\ref{1csums}).
 \end{remark}
 \begin{corollary}
For   $ |axt|<1$, we have:  
\begin{align}
\label{lums}
\sum_{n=0}^\infty \psi_n^{(a)}(x,y|q)(\lambda;q)_n \frac{(-1)^nq^{({}^n_2)}\,t^n}{(q;q)_n} 
= \frac{( \lambda, xt,yt; q)_\infty }{(axt;q)_\infty}  \;{}_{3}\Phi_2\left[\begin{array}{r}axt,0,0;
 \\\\
xt,yt; 
 \end{array} 
q; \lambda\right]
\end{align}
and
\begin{align}
\label{dd1sums}
\sum_{n=0}^\infty \psi_n^{(a)}(x|q)(\lambda;q)_n \frac{(-1)^nq^{({}^n_2)}\,t^n}{(q;q)_n} 
= \frac{( \lambda,xt, t; q)_\infty }{(axt;q)_\infty} \;{}_{3}\Phi_2\left[\begin{array}{r}axt,0,0;
 \\\\
xt,t ;
 \end{array}
q;\lambda\right].
\end{align}
 \end{corollary}
 \begin{remark}
Setting $\nu=1,\,a=0,\, y = ax$ and $z=y$ in (\ref{1sums}), we get (\ref{lums}). For   $\nu=1,\,a=0,\, y = ax$ and $z=1$,   (\ref{1sums}) reduces to   (\ref{dd1sums}). For  $y=1$, (\ref{lums}) reduces to  (\ref{dd1sums}).
 \end{remark}
 Before we prove  the Theorem  \ref{TA1}, the following Lemma is necessary.
 \begin{lemma}\cite[Eq. (III.1)]{GasparRahman90} For $\{|c|,|z|,|b|\}<1$, we have:
 \be
 \label{male}
 \;{}_{2}\Phi_1\left[\begin{array}{r}a,b; 
 \\\\
c;
 \end{array}
q;z\right]=\frac{(b,az;q)_\infty}{(c,z;q)_\infty}\;{}_{2}\Phi_1\left[\begin{array}{r}c/b,z ;
 \\\\
az;
 \end{array}
q; b\right].
\ee
 \end{lemma}
\begin{proof}[Proof of Theorems \ref{TA1}] 
Denoting the right-hand side of equation (\ref{1sums}) by $H'(a,t,x,y,z)$, we have:
\be 
\label{ass}
H'(a,t,x,y,z)=(\mu/\nu;q)_\infty\sum_{n=0}^\infty  \frac{(x\nu tq^n;q)_\infty\,(\mu/\nu)^n   }{(y\nu tq^n;q)_\infty(q;q)_n}\, \;{}_{1}\Phi_1\left[\begin{array}{r}a ; 
 \\\\
0;
 \end{array} 
q; z\nu tq^n\right].
\ee 
Because equation (\ref{ass}) satisfies (\ref{asaaaa}), by (\ref{shha}), we have:
\bea
\label{ss}
H'(a,t,x,y,z)&=&\widetilde{L}(a,z; \theta_{xy})   \left\{
H'(a,t,x,y,0)\right\} \cr
&=&\widetilde{L}(a,z; \theta_{xy})\Bigg\{ (\mu/\nu; q)_\infty\sum_{n=0}^\infty  \frac{(x\nu tq^n;q)_\infty\,(\mu/\nu)^n   }{(y\nu tq^n;q)_\infty(q;q)_n}\Bigg\}
\cr
&=& \widetilde{L}(a,z; \theta_{xy})\left\{  \frac{(\mu/\nu, x\nu t; q)_\infty }{(y\nu t;q)_\infty}\;{}_{2}\Phi_1\left[\begin{array}{r}y\nu t,0 ;\\
 \\
x\nu t;
 \end{array} 
q;\frac{\mu}{\nu}\right]\right\}.\nonumber 
\eea
By using   (\ref{male}) and (\ref{Srivas}), the last relation becomes
\bea
H'(a,t,x,y,z)
&=&\widetilde{L}(a,z; \theta_{xy})\Bigg\{\sum_{n=0}^\infty  
p_n (y,x)\frac{p_n(\nu,\mu)\,t^n}{(q;q)_n}\Bigg\} 
\cr
&=&\sum_{n=0}^\infty  
\widetilde{L}(a,z; \theta_{xy})\Bigg\{(-1)^nq^{-({}^n_2)} p_n (y,x)\Bigg\} p_n(\nu,\mu)\frac{(-1)^nq^{({}^n_2)}\,t^n}{(q;q)_n} \,\mbox{ by (\ref{conds})}
 \cr 
&=&\sum_{n=0}^\infty \Psi_n^{(a)}(x,y,z|q)p_n(\nu,\mu) \frac{(-1)^nq^{({}^n_2)}\,t^n}{(q;q)_n}, \nonumber
\eea
which is the left-hand side of  (\ref{1sums}).  This achieves  the proof.
 \end{proof}

\medskip

\section*{Acknowledgments}
The author would like to thank the referees and editors for their many valuable comments and suggestions. This work was supported by the Zhejiang Provincial Natural Science Foundation of China (No.~LY21A010019).

\end{CJK*}
\end{document}